\documentclass{article}

\usepackage{graphics,psfrag,epsfig,amsbsy,amsfonts,enumerate,amsmath,amsthm} 
 
\def\BV{\mathinner{\bf BV}}
\def\L#1{\mathinner{{{\bf L}^{ #1}}}}
\def\C#1{\mathinner{{{\bf C}^{\bf #1}}}}

\def\reali{{\mathbb R}}
\def\interi{{\mathbb Z}}

\def\tv{\mathop{\mbox {\rm Tot.Var.}}}

\def\norma#1{{\left\| #1 \right\|}}

\def \dd #1 #2 {{\frac {\partial #1} {\partial #2} }}

\def\ff{\phi}

\def \1 {{\bf 1}}
\def \ss {{\boldsymbol {\sigma} }}

\def \al {\alpha}

\def \la {\lambda}
\def \ve {\varepsilon}

\def \s  {\sigma}

\def \Om {\Omega}

\def \ka{{k_\alpha}}
\def \cD{{\cal D}}

\def \NP{{\cal NP}}

\def \J{{\cal J}}

\def \Z{{\cal Z}}

\def \O{{\cal O}}
\def \R{{\cal R}}
\def \S{{\cal S}}

\long\def\rem#1{\null}

\def \non {\nonumber}
\def \pn{\par \noindent}

\def\d{{\delta}}

\def\om{\omega}

\def \eps{\epsilon}

\def\sqr#1#2{\vbox{
   \hrule height .#2pt 
   \hbox{\vrule width .#2pt height #1pt \kern #1pt 
      \vrule width .#2pt}
   \hrule height .#2pt }}
\def\square{\sqr42}

\def\sddots{\mathinner{\mkern1mu\lower7pt\vbox{\kern7pt\hbox{.}}\mkern2mu
\lower4pt\hbox{.}\mkern2mu\lower1pt\hbox{.}\mkern1mu}}

\def\buildrel#1\over#2{\mathrel{\mathop{\kern0pt #2}\limits^{#1}}}

\def\be{\begin{eqnarray}}
\def\ee{\end{eqnarray}}
\def\ben{\begin{enumerate}}
\def\een{\end{enumerate}}
\def\ba{\begin{array}}
\def\ea{\end{array}}

\newtheorem{Theorem}{Theorem}
\newtheorem{Proposition}{Proposition}
\newtheorem{Lemma}{Lemma}
\theoremstyle{definition}
\newtheorem{Remark}{Remark}
\newtheorem{Definition}{Definition}

\newcommand{\dref}[1]{(\ref{#1})}

\begin{document}

\title{Balance laws with integrable unbounded sources 
\thanks{This work has been supported by the SPP1253 priority program of the DFG and by the DAAD program D/06/19582.}
} 

\author{Graziano Guerra\footnote{Universit\`{a} di Milano--Bicocca, Italy, graziano.guerra@unimib.it} \and Francesca Marcellini\footnote{Universit\`{a} di Milano--Bicocca, Italy, f.marcellini@campus.unimib.it} \and Veronika Schleper\footnote{TU Kaiserslautern, Germany, sachers@mathematik.uni-kl.de}}

\maketitle 
\begin{abstract} 
We consider the Cauchy problem for a $n\times n$ strictly hyperbolic system of 
balance laws
$$
\left\{\begin{array}{c}
u_t+f(u)_x=g(x,u), \qquad x \in \reali, t>0\\
u(0,.)=u_o \in \L1 \cap \BV(\reali; \reali^n), \\
| \la_i(u)| \geq c > 0 \mbox{ for all } i\in \{1,\ldots,n\}, \\
\norma{g(x,\cdot)}_{\mathbf{C}^2}\leq \tilde M(x) \in \L1, \\
\end{array}\right.
$$
each characteristic field being genuinely nonlinear or linearly degenerate. 
Assuming that the $\mathbf{L}^1$ norm of $\|g(x,\cdot)\|_{\mathbf{C}^1}$ and $\|u_o\|_{\BV(\reali)}$ are small
enough, we prove the existence and uniqueness of global entropy solutions of 
bounded total variation extending the result in \cite{AmaGoGue} to unbounded (in $\L\infty$) sources. 
Furthermore, we apply this result to the fluid flow in a pipe with discontinuous cross sectional area, showing existence and uniqueness of the underlying semigroup.

  \medskip

  \noindent\textit{2000~Mathematics Subject Classification:} 35L65,
   35L45, 35L60.

  \medskip


\end{abstract} 

\section{Introduction}

The recent literature offers several results on the properties
of gas flows on networks. For instance,
in~\cite{ColomboGaravello2007,ColomboGuerraHertySachers,ColomboHertySachers,ColomboMauri} the well posedness is established for the gas flow at a junction of $n$ pipes with constant diameters. The equations governing the gas flow in a pipe with a smooth varying cross section $a(x)$ are given by (see for instance \cite{Liu79}):

\begin{displaymath}
\left\{ 
\begin{array}{l}
   \frac{\partial\rho}{\partial t}+ \frac{\partial q}{\partial 
       x}=-\frac{a'(x)}{a(x)}q\\
   \frac{\partial q}{\partial t}+ \frac{\partial (\frac{q^2}{\rho}+p)}{\partial 
       x}=-\frac{a'(x)}{a(x)}\frac{q^2}{\rho}\\ 
   \frac{\partial e}{\partial t}+ \frac{\partial \left(\frac{q}{\rho}(e+p)\right)}{\partial 
       x}=-\frac{a'(x)}{a(x)}\left(\frac{q}{\rho}(e+p)\right).\\ 
 \end{array}
\right.
\end{displaymath}
The well posedness of this system is covered in \cite{AmaGoGue} where 
an attractive unified approach to the
existence and uniqueness theory for quasilinear strictly hyperbolic
systems of balance laws is proposed. The case of discontinuous cross sections is considered in the literature inserting a junction with suitable coupling conditions at the junction, see for example \cite{ColomboGaravello2007,ColomboGuerraHertySachers,lefgoa}. 
One way to obtain coupling conditions at the point of discontinuity of the cross section $a$ is to take the limit of a sequence of Lipschitz continuous cross sections $a^\varepsilon$
converging to $a$ in $\L1$
(for a different approach see for instance \cite{ColMar}).
Unfortunately the results in \cite{AmaGoGue} require $\L\infty$ bounds on the source term and well posedness is proved on a domain depending on this $\L\infty$ bound. Since in the previous equations the source term contains the derivative of the cross sectional area one cannot hope to take the limit $a^\varepsilon\to a$. Indeed when $a$ is discontinuous, the $\L\infty$ norm of $(a^\varepsilon)'$ goes to infinity.
Therefore the purpose of this paper is to establish the result in \cite{AmaGoGue} without requiring the $\L\infty$ bound.
More precisely, we consider
the Cauchy problem for the following $n\times
n$ system of equations
\begin{equation}
   u_t+f(u)_x=g(x,u), \qquad x\in \reali,\ t>0, \label{integrable:eq:bal}
\end{equation}
endowed with a (suitably small) initial data
\begin{equation}
        u(0,x)=u_o(x), \qquad x\in \reali.  \label{init:data}
\end{equation}
 belonging to $\L1 \cap
\BV \left(\reali;\reali^n\right)$, the
space of integrable functions with bounded total variation (Tot.Var.)
in the sense of \cite{vol}.
Here $u(t,x)\in \reali^n$ is the vector of unknowns, $f:\Om \to \reali^n$
denotes the fluxes, {\it i.e.} a smooth function defined on
$\Om$ which is an open neighborhood of the origin in $\reali^n$. The
system (\ref{integrable:eq:bal}) is supposed to be strictly
hyperbolic, with each characteristic field either genuinely nonlinear
or linearly degenerate in the sense of Lax \cite{Lax}. Concerning the
source term $g$, we assume that it satisfies the following
Caratheodory--type conditions:
\smallskip
\begin{itemize}
\item[$(P_1)$]
$g: \reali\times\Om \to \reali^n$ is measurable with respect to (w.r.t.)
$x$, for any $u\in\Om$, and is $\C2$ w.r.t. $u$, for any $x\in\reali$;
\item[$(P_2)$] there exists a $\L1$ function $\tilde M(x)$ such that $\norma{g(x,\cdot)}_{\C2}\le \tilde M(x)$;
\label{hyp_on_g}
\item[$(P_3)$]
 there exists a function  $\omega \in \L1 (\reali)$
such that $\|g(x,\cdot)\|_{\C1}\leq \om(x)$.
\end{itemize}
\begin{Remark}
 Note that the $\L1$ norm of $\tilde M(x)$ does not have to be small but only bounded differently from $\omega(x)$ whose norm has to be small (see Theorem \ref{th:introduction} below). Furthermore condition $(P_2)$ replaces the $\L\infty$ bound of the $\C2$ norm of $g$ in \cite{AmaGoGue}. Finally observe that we do not require any $\L\infty$ bound on $\omega$. On the other hand we will need the following observation: if we define 
\begin{equation}\label{epsilonh}
  \tilde\varepsilon_h=\sup_{x\in\reali}\int_0^h\omega(x+s)\;ds,
\end{equation}
by absolute continuity one has $\tilde \varepsilon_h\to 0$ as $h\to 0$.

\end{Remark}
\smallskip
Moreover, we assume that a {\sl non-resonance} condition holds, that
is the characteristic speeds of the system \dref{integrable:eq:bal}
are bounded away from zero:
\be
| \la_i(u)| \geq c > 0, \qquad \forall~u \in \Om,\ i\in \{1,\ldots,n\}.
\label{hyp:non_resonance}
\ee

The following theorem states the well posedness of (\ref{integrable:eq:bal}) in the above defined setting.

\begin{Theorem}\label{th:introduction}
Assume $(P_1)$--$(P_3)$ and \dref{hyp:non_resonance}. If the 
norm of $\omega$ in $\L1(\reali)$ is sufficiently small, 
there exist a constant $L>0$,
a closed domain $\cD$ of integrable functions with small total
variation and a unique semigroup $P:[0,+\infty)\times
\cD\rightarrow\cD$ satisfying
\begin{enumerate}[i)]
\item
$P_0u=u,\quad P_{t+s}u=P_t \circ P_su$ for all $u,v\in\cD$ and $t,s\ge 0$;
\item
$\|P_su-P_tv\|_{\L1 (\reali)}\le L\Big(|s-t|+\|u-v\|_{\L1
(\reali)}\Big)$ for all $u,v\in\cD$ and $t,s\ge 0$;
\item
for all $u_o\in\cD$ the function $u(t,\cdot)=P_tu_o$ is a weak entropy
solution of the Cauchy problem
\dref{integrable:eq:bal}--\dref{init:data} and satisfies the integral
estimates \dref{eq:first_int_cond}, \dref{eq:second_int_cond}.
\end{enumerate}

Conversely let $u:[0,T]\rightarrow\cD$ be Lipschitz continuous as a map
with values in $\L1(\reali,\reali^n)$ and assume that $u(t,x)$
satisfies the integral conditions \dref{eq:first_int_cond},
\dref{eq:second_int_cond}. Then $u(t,\cdot)$ coincides with a
trajectory of the semigroup $P$.
\end{Theorem}
The proof of this theorem is postponed to sections \ref{sec:h-riemann_solver} and \ref{sec:uniqueness}, where existence and uniqueness are proved. Before these technical details, we state the application of the above result to gas flow in section \ref{application}. Here we apply Theorem \ref{th:introduction} to establish the existence and uniqueness of the semigroup related to pipes with discontinuous cross sections. Furthermore, we show that our approach yields the same semigroup as the approach followed in \cite{ColomboHertySachers} in the special case of two connected pipes. The technical details of section \ref{application} can be found at the end of the paper in section \ref{section5}.

\section{Application to gas dynamics}
\label{application}

\begin{figure}[h]
\centering
 \includegraphics[width=0.4\textwidth]{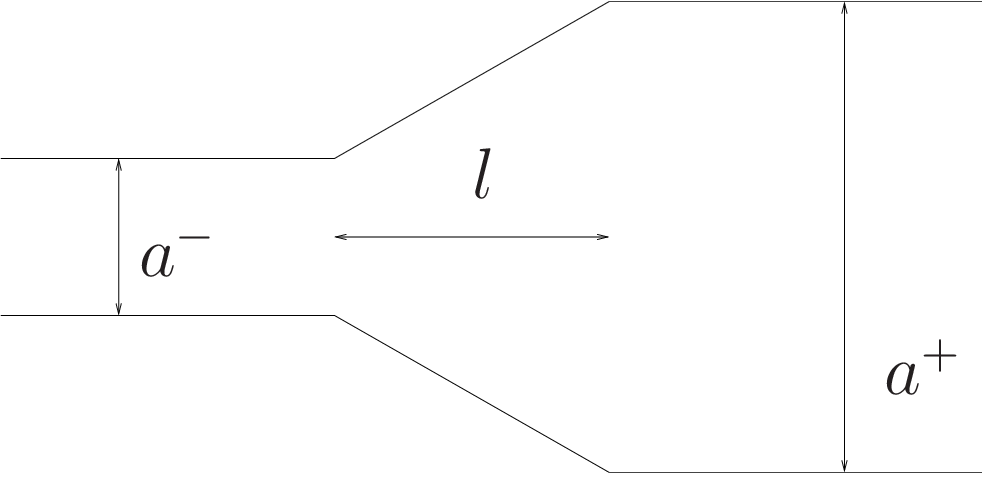}
 \hspace*{0.09\textwidth}
 \includegraphics[width=0.4\textwidth]{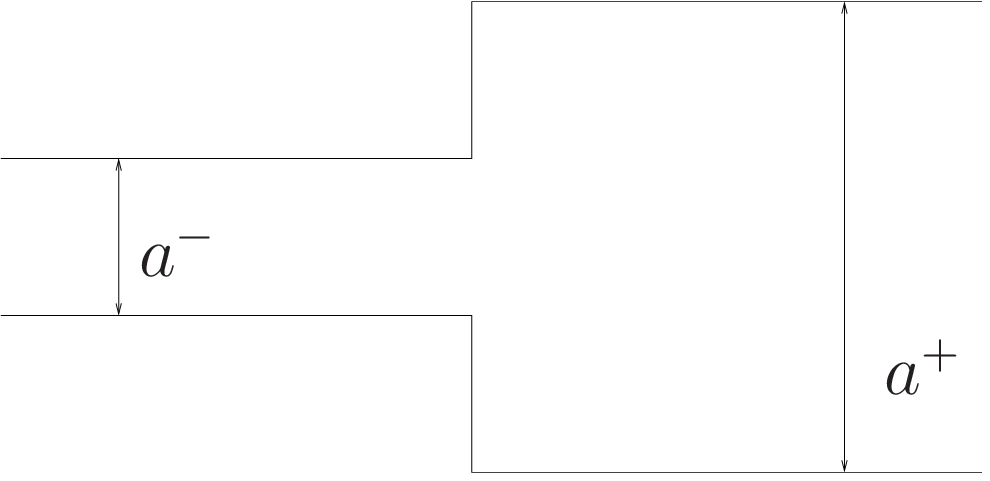}
\caption{Illustration of approximated and discontinuous cross-sectional area}
\label{fig:1}
\end{figure}

Theorem~\ref{th:introduction} provides an existence and uniqueness result for pipes with Lipschitz continuous cross section where the equations governing the gas flow are given by
\begin{equation}\label{dynamics}
\left\{ 
\begin{array}{l}
   \frac{\partial\rho}{\partial t}+ \frac{\partial q}{\partial 
       x}=-\frac{a'(x)}{a(x)}q\\
   \frac{\partial q}{\partial t}+ \frac{\partial (\frac{q^2}{\rho}+p)}{\partial 
       x}=-\frac{a'(x)}{a(x)}\frac{q^2}{\rho}\\ 
   \frac{\partial e}{\partial t}+ \frac{\partial \left(\frac{q}{\rho}(e+p)\right)}{\partial 
       x}=-\frac{a'(x)}{a(x)}\left(\frac{q}{\rho}(e+p)\right).
 \end{array}
\right. 
\end{equation}
Here, as usual, $\rho$ denotes the mass density, $q$ the linear momentum, $e$ is the energy density, $a$ is the area of the cross section of the pipe and $p$ is 
the pressure which is related to the conserved quantities $(\rho,q,e)$ by the equations of state.
In most situations, when two pipes of different size have to be connected, the length $l$ of the adaptor is small compared to the length of the pipes. Therefore it is convenient to model these connections as pipes with a jump in the cross sectional area. These discontinuous cross sections however do not fulfill the requirements of Theorem~\ref{th:introduction}. Nevertheless, we can use this Theorem to derive the existence of solutions to the discontinuous problem by a limit procedure.
To this end, we approximate the discontinuous function
\begin{equation}
 a(x) = \left\{\begin{array}{ll} a^-, & x<0\\
                a^+, & x>0
               \end{array}\right.
\end{equation}
by a sequence $a_l\in C^{0,1}({\mathbb R},{\mathbb R}^+)$ with the following properties
\begin{equation}
 a_l(x) = \left\{\begin{array}{ll} a^-, & x<-\frac{l}{2} \\
                  		   \varphi_l(x), & x\in\left[-\frac{l}2, \frac{l}2\right]\\
				   a^+, & x>\frac{l}2
                 \end{array}\right.
\end{equation}
where $\varphi_l$ is any smooth monotone function which connects the two strictly positive constants $a^-$, $a^+$.
One possible choice of the approximations $a_l$ as well as the discontinuous pipe with cross section $a$ are shown in figure~\ref{fig:1}.

With the help of Theorem~\ref{th:introduction} and the techniques used in its proof, we are now able to derive the following Theorem (see also \cite{ColMar} for a similar result obtained with different methods).

\begin{Theorem}\label{Theorem2}
 If  $\|a_l^\prime\|_{\L1}=|a^+-a^-|$ is sufficiently small, the semigroups $P^l$ related with the smooth section $a_l$ converge to a unique semigroup $P$. 

The limit semigroup satisfies and is uniquely identified by the integral estimates \dref{eq:first_int_cond},
\dref{eq:second_int_cond} with $U^\sharp$ substituted by $\bar U^\sharp$ (see Section \ref{section5}) for the point $\xi=0$.
More precisely let $u:[0,T]\rightarrow\cD$ be Lipschitz continuous as a map
with values in $\L1(\reali,\reali^n)$ and assume that $u(t,x)$
satisfies the integral conditions \dref{eq:first_int_cond},
\dref{eq:second_int_cond} with $U^\sharp$ substituted by $\bar U^\sharp$ for the point $\xi=0$. Then $u(t,\cdot)$ coincides with a
trajectory of the semigroup $P$.
\end{Theorem}

Observe that the same Theorem holds for the $2\times 2$ isentropic 
system (see Section \ref{section5}) 
\begin{equation}\label{dynamic}
\left\{ 
\begin{array}{l}
   \frac{\partial\rho}{\partial t}+ \frac{\partial q}{\partial 
       x}=-\frac{a'(x)}{a(x)}q\\
   \frac{\partial q}{\partial t}+ \frac{\partial (\frac{q^2}{\rho}+p)}{\partial 
       x}=-\frac{a'(x)}{a(x)}\frac{q^2}{\rho}. 
\end{array}
\right. 
\end{equation}
In \cite{ColomboHertySachers} $2\times 2$ homogeneous conservation laws at a junction are considered for given admissible junction conditions. The situation of a junctions with only two pipes with different cross sections can be modeled by our limit procedure or as in  \cite{ColomboHertySachers} with a suitable junction condition. If we define the function $\Psi$ which describes the junction conditions as 
\begin{equation}\label{eq:psieq}
 \Psi\left(\rho_1,q_1,\rho_2,q_2\right)=
(\rho_2,q_2)-\Phi(a^+-a^-,\rho_1,q_1)
\end{equation}
then it fulfills the determinant condition  in \cite[Proposition 2.2]{ColomboHertySachers} since it satisfies Lemma \ref{lemma:Riemann_prob0}. Here $\Phi(a,u)$ is the solution of the ordinary differential equation (\ref{ODE}) in Section \ref{section5}.
With these junction conditions one can show that the semigroup obtained in \cite{ColomboHertySachers} satisfies the same integral estimate (see the following proposition) as our limit semigroup hence they coincide.

\begin{Proposition}\label{prop1}
 The semigroup defined in \cite{ColomboHertySachers} with the junction condition given by 
 (\ref{dynamics}) satisfies the integral estimates \dref{eq:first_int_cond},
\dref{eq:second_int_cond} with $U^\sharp$ substituted by $\bar U^\sharp$ for the point $\xi=0$.
\end{Proposition}

The proof is postponed to Section \ref{section5}.

\begin{Remark}
 Note that Proposition 1 justifies the coupling condition (\ref{eq:psieq}) as well as the condition 
used in  \cite{lefgoa} to study the Riemann problem for the gas flow through a nozzle. 
\end{Remark}

\section{Existence of BV entropy solutions}
\label{sec:h-riemann_solver}

Throughout the next two sections, we follow the structure of \cite{AmaGoGue}. We recall some definitions and notations in there, and also the results which do not depend on the $\L\infty$ boundedness of the source term. We will prove only the results which in \cite{AmaGoGue} do depend on the $\L\infty$ bound using our weaker hypotheses.

\subsection{The non homogeneous Riemann-Solver}
Consider the stationary equations associated to
\dref{integrable:eq:bal}, namely the system of ordinary differential
equations:
\be
   f(v(x))_x=g(x,v(x)). \label{eq:ode}
\ee
For any $x_o\in\reali$, $v\in\Om$, consider the initial
data
\be
v(x_o) = v.     \label{eq:init-data}
\ee
As in \cite{AmaGoGue}, we introduce a suitable approximation of the
solutions to \dref{eq:ode}, \dref{eq:init-data}. Thanks to
\dref{hyp:non_resonance}, the map $u \mapsto f(u)$ is invertible
inside some neighborhood of the origin; in this neighborhood,
for small $h>0$, we can define
\be
\Phi_h(x_o,u)\,\dot=\,
f^{-1}\left[f(u) + \int_0^h g\left(x_o+s,u\right)\,ds \right].
\label{def:Phi}
\ee
This map gives an approximation of the flow of \dref{eq:ode} in the
sense that
\be
f\left(\Phi_h(x_o,u)\right)-f(u)=\int_0^hg\left(x_o+s,u\right)\,ds.
\label{approx_ode}
\ee

Throughout the paper we will use the Landau notation
$\O(1)$ to indicate any function whose absolute value remains
uniformly bounded, the bound depending only on $f$ and $\|\tilde M\|_{\L1}$.

\begin{Lemma}
\label{lemma0}
 The function $\Phi_h(x_o,u)$ defined in (\ref{def:Phi}) satisfies the following 
 uniform (with respect to $x_o\in \reali$ and to $u$ in a suitable neighborhood of 
 the origin) estimates.
 \begin{equation}
 \begin{split}
  &\|\Phi_h(x_o,\cdot)\|_{\C2}\le\O(1),\qquad \lim_{h\to 0}\sup_{x_o\in \reali}
  |\Phi_h(x_o,u)-u|=0,\\ 
  &\lim_{h\to 0}\|Id - D_u\Phi_h(x_o,u)\|=0
 \end{split}
 \end{equation}
\end{Lemma}

\begin{proof}
 The Lipschitz continuity of $f^{-1}$ and (\ref{epsilonh}) imply
 \begin{displaymath}
  \begin{split}
   \left|\Phi_h(x_o,u)-u\right|&=
   \left|\Phi_h(x_o,u)-f^{-1}\left(f(u)\right)\right|\le
   \O(1)\left|\int_0^{h}g(x_o+s,u)\;ds\right|\\
   &\le\O(1)\left|\int_0^{h}\omega(x_o+s)\;ds\right|
   \le\O(1)\tilde\varepsilon_h\xrightarrow{h\to 0}0.
  \end{split}
 \end{displaymath}
 Next we compute 
 \begin{displaymath}
  \begin{split}
   D_u\Phi_h(x_o,u)=&Df^{-1}\left[f(u)+\int_0^{h}g(x_o+s,u)\;ds\right]\\
   &\qquad\cdot \left(Df(u)+\int_0^{h}D_ug(x_o+s,u)\;ds\right)
  \end{split}
 \end{displaymath}
 which together to the identity $u=f^{-1}\left(f(u)\right)$ implies
 \begin{displaymath}
 \begin{split}
  \left\|D_u\Phi_h(x_o,u)-Id\right\|&=\left\|D_u\Phi_h(x_o,u)-
  Df^{-1}\left(f(u\right)\right)\|\\
  &\le \left\|Df^{-1}\left[f(u)+\int_0^{h}g(x_o+s,u)\;ds\right]-Df^{-1}\left(f(u)\right)\right\|\\
  &\qquad\qquad\cdot\left(\|Df(u)\|+\int_0^{h}\|D_ug(x_o+s,u)\|\;ds\right)\\
  &\qquad +\left\|Df^{-1}\left(f(u)\right)\right\|\cdot\int_0^{h}\|D_ug(x_o+s,u)\|\;ds\\
  &\le \O(1)\tilde\varepsilon_h\xrightarrow{h\to 0}0.
  \end{split}
 \end{displaymath}
Finally, denoting with $D_i$ the partial derivative with respect to the $i$ component of the state vector and by $\Phi_{h,\ell}$ the $\ell$ component of the vector $\Phi_h$, we derive
\begin{displaymath}
 \begin{split}
  D_iD_j\Phi_{h,l}(x_o,u)&=\sum_{k,k'}\Bigg(
  D_kD_{k'}f_\ell^{-1}\left(f(u)+\int_0^hg(x_o+s,u)\;ds\right)\\  
  &\qquad\quad\cdot\left(D_if_k(u)+\int_0^hD_ig_k(x_o+s,u)\;ds\right)\\
  &\qquad\quad\cdot\left(D_jf_{k'}(u)+\int_0^hD_jg_{k'}(x_o+s,u)\;ds\right)\Bigg)\\
  &\qquad +\sum_{k}D_kf^{-1}_\ell\left(f(u)+\int_0^hg(x_o+s,u)\;ds\right)\\
  &\qquad\quad\cdot\left(D_jD_if_k(u)+\int_0^hD_jD_ig_k(x_o+s,u)\;ds\right)
 \end{split}
\end{displaymath}
so that
\begin{displaymath}
 \begin{split}
  \left\|D^2\Phi_{h}(x_o,u)\right\|&\le\O(1)\left(1+\int_0^h\tilde M(x_o+s)\;ds\right)\le\O(1)\left(1+\|\tilde M\|_{\L1}\right)\le\O(1).
 \end{split}
\end{displaymath}

\end{proof}

For any $x_o\in\reali$ we
consider the system \dref{integrable:eq:bal},
endowed with a Riemann initial datum:
\be
u(0,x)= \left\{
\begin{array}{ll}
u_\ell  & {\rm if }\  x<x_o\\
u_r     & {\rm if }\  x>x_o.
\end{array}
\right.
\label{eq:riemann_data}
\ee
If the two states $u_\ell$,
$u_r$ are sufficiently close, let $\Psi$ be the unique entropic homogeneous
Riemann solver given by the map
\be
u_r = \Psi(\ss)(u_\ell) =
\psi_n(\s_n)\circ\ldots\circ\psi_1(\s_1)(u_\ell), \non
\ee
where $\ss=(\s_1,\ldots,\s_n)$ denotes the (signed) wave strengths vector
in $\reali^{n}$, \cite{Lax}. Here $\psi_j$, $j=1,\ldots,n$ is the
shock--rarefaction curve of the $j^{th}$ family, parametrized as
in \cite{Bdue} and related to the homogeneous system of conservation laws
\be
        u_t+f(u)_x=0. \label{eq:homo}
\ee
Observe that, due to \dref{hyp:non_resonance}, all the simple waves
appearing in the solution of \dref{eq:homo}, \dref{eq:riemann_data}
propagate with {\it non-zero} speed.

\smallskip
To take into account the effects of the source term, we consider a
stationary discontinuity across the line $x=x_o$, that is, a wave whose
speed is equal to 0, the so called \underline{zero-wave}.
Now, given $h>0$, we say that the particular Riemann solution:
\be u(t,x)= \left\{
\begin{array}{ll}
u_\ell  & {\rm if }\  x<x_o\\
u_r     & {\rm if }\  x>x_o.
\end{array}
\right. \qquad \forall~t\geq 0
\label{zero-wave}
\ee
is admissible if and only if $u_r = \Phi_h (x_o,u_\ell)$, where
$\Phi_h$ is the map defined in \dref{def:Phi}. Roughly speaking, we
require $u_\ell$, $u_r$ to be (approximately) connected by a solution
of the stationary equations \dref{eq:ode}.

\begin{Definition}\label{def:h_R_solver}
Given $h>0$ suitably small, $x_o\in\reali$, we say that $u(t,x)$ is a
{$h$--Rie\-mann sol\-ver} for
\dref{integrable:eq:bal}, \dref{hyp:non_resonance}, \dref{eq:riemann_data},
if the following conditions hold
\begin{itemize}
\item[(a)] there exist two states $u^-$, $u^+$ which satisfy
$u^+ = \Phi_h(x_o,u^-)$;
\item[(b)]
on the set $\{t\geq 0,\ x<x_o\}$, $u(t,x)$ coincides with the
solution to the homogeneous Riemann Problem \dref{eq:homo}
with initial values $u_\ell$, $u^-$ and,
on the set $\{t\geq 0,\ x>x_o\}$, with the solution to the homogeneous
Riemann Problem with initial values $u^+$, $u_r$;
\item[(c)] the Riemann Problem between $u_\ell$ and $u^-$ is solved only
by waves with negative speed (i.e. of the families $1,\ldots,p$);
\item[(d)]  the Riemann Problem between $u^+$ and $u_r$ is solved only
by waves with positive speed (i.e. of the families $p+1,\ldots,n$).
\end{itemize}
\end{Definition}
\pn

\begin{figure}[ht]
\centering
\epsfig{figure=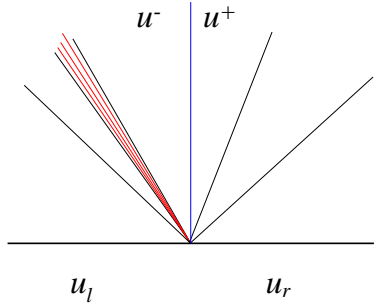,height=4cm}
\caption{Wave structure in an $h$--Riemann solver.}
\label{HRS}
\end{figure}

\begin{Lemma}
\label{initial_estimates}
Let $x_o\in\reali$ and $u,\ u_1,\ u_2$ be three states in a suitable
neighborhood of the origin. For $h$ suitably small, one has
\be
\label{eq:useful_est}
\left|\Phi_h(x_o,u)-u\right| &=& \O(1)\int_0^h\om(x_o+s)\;ds, \\
\left|\Phi_h(x_o,u_2)-\Phi_h(x_o,u_1)-(u_2 - u_1) \right|&=& 
\O(1)|u_2-u_1|\int_{x_o}^{x_o+h}\om(s)ds.
\label{eq:useful_est2}
\ee
\end{Lemma}

\begin{Lemma}
\label{lemma:Riemann_prob0}
For any $M>0$ there exist $\d_1^\prime,\ h_1^\prime>0$, depending only on
$M$ and
the homogeneous system \dref{eq:homo}, such that the following holds.
For all maps $\ff\in\C2 (\reali^n,\reali^n)$ satisfying
$$
\|\ff\|_{\C2}\le M,\qquad
        |\ff(u)-u|\le h_1^\prime,\qquad
                \| I-D\ff(u)\|\le h_1^\prime
$$
and for all $u_\ell\in B(0,\d_1^\prime)$, $u_r\in B(\ff(0),\d_1^\prime)$
there exist
$n+1$ states $w_0,\ldots,w_{n+1}$ and $n$ wave sizes
$\s_1,\ldots,\s_n$, depending smoothly on $u_\ell,\ u_r$, such that with
previous notations:
\begin{enumerate}[i)]
\item
$w_0=u_\ell,\ w_{n+1}=u_r$;
\item
$w_i=\Psi_i(\s_i)(w_{i-1}),\qquad i=1,\ldots,p$;
\item
$w_{p+1}=\ff(w_p)$;
\item
$w_{i+1}=\Psi_i(\s_{i})(w_{i}),\qquad i=p+1,\ldots,n$.
\end{enumerate}
\end{Lemma}

The next lemma establishes existence and uniqueness for the
$h$--Riemann solvers (see Fig.\ref{HRS}).

\begin{Lemma}
\label{lemma:Riemann_prob}
There exist $\d_1,\ h_1>0$ such that the following holds: for any $x_o\in
\reali$, $h\in [0,h_1]$, $u_\ell$, $u_r \in B(0,\d_1)$, there exists
a unique $h$--Riemann solver in the sense of Definition~\ref{def:h_R_solver}.
\end{Lemma}

\begin{proof}
By Lemma
\ref{lemma0} if $h_1>0$ is chosen
sufficiently small then for any $h\in[0,h_1]$, $x_o\in\reali$ the map
$u\mapsto
\Phi_h(x_o,u)$ meets the hypotheses of Lemma
\ref{lemma:Riemann_prob0}.
Finally taking $h_1$ eventually smaller we can obtain that there exists
$\d_1>0$ such that $B(0,\d_1)\subset B(0,\d_1^\prime)\cap
B(\Phi_h(x_o,0),
\d_1^\prime)$, for any $h\in [0,h_1]$.
\end{proof}

\smallskip
In the sequel, $E$ stands for the implicit function given by Lemmas
\ref{lemma:Riemann_prob0} and \ref{lemma:Riemann_prob}:
$$
\ss \dot=E(h,u_\ell,u_r;x_o),
$$
which plays the role of a wave--size vector. We recall that, by Lemma \ref{lemma:Riemann_prob0}, $E$ is a
$\C2$ function with respect to the variables $u_\ell,u_r$ and its
$\C2$ norm is bounded by a constant independent of $h$ and $x_o$.

In contrast with the homogeneous case, the wave--size
$\ss$ in the $h$--Riemann solver is not equivalent to the jump size
$|u_\ell-u_r|$; an additional term appears coming from the ``Dirac source
term" (see the special case $u_\ell=u_r$).

\begin{Lemma}\label{lemma:equiv_tot_var}
Let $\d_1,\ h_1$ be the constants in Lemma
\ref{lemma:Riemann_prob}. For $u_\ell,u_r\in B(0,\d_1)$, $h\in [0,h_1]$,
set $\ss=E[h,u_\ell,u_r;x_o]$. Then it holds:
\be
\begin{array}{rcl}
|u_\ell-u_r|&=&\O(1)\left(|\ss|+\displaystyle{\int_0^h\om(x_o+s)\;ds}
\right),\\
|\ss|&=&\O(1)\left(|u_\ell-u_r|+
\displaystyle{\int_0^h\om(x_o+s)\;ds}\right).
\end{array}
\label{eq:equiv_tot_var}
\ee
\end{Lemma}

\subsection{Existence of a Lipschitz semigroup of $\BV$ entropy solutions}
\label{sec:existence_entropy_solutions}
Note that as shown in \cite{AmaGoGue} we can identify the sizes of the zero waves with the quantity
\be
   \s = \int_{0}^{h} \om(jh+s)\,ds.
\label{size_0_wave}
\ee
With this definition all the Glimm interaction estimates continue to hold with constants that depend only on $f$ and on $\|\tilde M\|_{\L1}$, therefore all the wave front tracking algorithm can be carried out obtaining the existence of $\varepsilon,h$-approximate solutions as defined below.

\begin{Definition}\label{decadix}
Given $\eps, h >0$, we say that a continuous map
$$
u^{\eps,h}:
[0,+\infty)\rightarrow \L1_{loc}\left(\reali,\reali^n\right)
$$
is an $\eps,h$--approximate solution of
\dref{integrable:eq:bal}--\dref{init:data} if the following holds:
\begin{itemize}

\item[--] As a function of two variables, $u^{\eps,h}$ is piecewise
constant with discontinuities occurring along finitely many straight
lines in the $x,t$ plane.  Only finitely many wave-front interactions
occur, each involving exactly two wave-fronts, and jumps can be of
four types: shocks (or contact discontinuities), rarefaction waves,
non-physical waves and zero-waves: $\J=\S \cup \R \cup \NP \cup \Z$.

\item[--] Along each shock (or contact discontinuity) $x_\al=x_\al(t)$,
$\al \in \S$, the values of $u^-=u^{\eps,h}(t,x_\al -)$ and
$u^+=u^{\eps,h}(t,x_\al +)$ are related by
$u^+=\psi_{\ka}(\sigma_\al)(u^-)$ for some $\ka \in \{1, ..., n\}$ and
some wave-strength $\sigma_\al$. If the $\ka^{th}$ family is genuinely
nonlinear, then the Lax entropy admissibility condition $\sigma_\al
<0$ also holds. Moreover, one has
$$
|\dot{x}_\al - \lambda_{\ka}(u^+,u^-)| \leq \eps
$$
where $\lambda_{\ka}(u^+,u^-)$ is the speed of the
shock front (or contact
discontinuity) pre\-scribed by the classical
Ran\-ki\-ne-Hu\-goniot conditions.

\item[--] Along each rarefaction front $x_\al =x_\al(t)$, $\al \in \R$,
one has $u^+=\psi_{\ka}(\sigma_\al)(u^-)$, $0 < \sigma_\al \leq
\eps$ for some genuinely nonlinear family $\ka$. Moreover, we have:
$|\dot{x}_\al - \lambda_{\ka}(u^+)| \leq \eps$.

\item[--] All non-physical fronts $x=x_\al(t)$, $\al \in \NP$ travel
at the same speed $\dot{x}_\al = \hat{\lambda} >
\sup_{u,i}|\lambda_i(u)|$. Their total strength remains uniformly
small, namely:
$$
\sum_{\al \in \NP}
|u^{\eps,h}(t,x_\al +)-u^{\eps,h}(t,x_\al -)| \leq \eps,\qquad
\forall ~t > 0.
$$

\item[--] The zero-waves are located at every point $x=jh$, $j\in
(-\frac{1}{h\eps},\frac{1}{h\eps}) \cap \interi$.\\ Along a
zero-wave located at $x_\alpha=j_\alpha h$, $\al\in\Z$,
the values $u^-=u^{\eps,h}(t,x_\al -)$ and
$u^+=u^{\eps,h}(t,x_\al +)$ satisfy $u^+=\Phi_h(x_\al,u^-)$ for all
$t>0$ except at the interaction points.

\item[--] The total variation in space $\tv u^{\eps,h}(t,\cdot)$
is uniformly bounded for all $t\ge 0$.
The total variation in time $\tv
\left\{u^{\eps,h}(\cdot,x);[0,+\infty)\right\}$
is uniformly bounded for $x\not= jh$, $j \in \interi$.
\end{itemize}

Finally, we require that $\|u^{\eps,h}(0,.)-u_o\|_{\L1(\reali)} \leq \eps$.
\end{Definition}

Keeping $h>0$ fixed, we are about to let first
$\eps$ tend to zero. Hence we shall drop the superscript $h$
for notational clarity.

\begin{Theorem}\label{first-convergence}
Let $u^{\eps}$ be a family of $\eps,h$--approximate solutions of
\dref{integrable:eq:bal}--\dref{init:data}. There exists a subsequence
$u^{\eps_i}$ converging as $i \rightarrow +\infty$ in ${\bf
L}^1_{loc}\left((0,+\infty)\times\reali\right)$ to a function $u$
which satisfies for any $\varphi\in {\bf C}^1_c
\left((0,+\infty)\times\reali\right)$:
\be
&&\int_0^\infty\int_\reali \left[u\varphi_t+f(u)\varphi_x\right]
dxdt \non \\
&&+\int_0^{\infty}\sum_{j\in\interi}\varphi(t,jh)\left(\int_0^h
g\left[jh+s,u(t,jh-)\right]ds\right)dt =0.\label{eq:hequation}
\ee
Moreover $\tv u(t,\cdot)$ is uniformly bounded and $u$ satisfies the
Lipschitz property
\begin{equation}
\int_\reali\left|u(t',x)-u(t'',x)\right|dx\le C'|t'-t''|,\qquad
t',t''\ge 0;
\label{eq:lip_in_t}
\end{equation}
\end{Theorem}

\smallskip Now we are in position to prove \cite[Theorem 4]{AmaGoGue} with our weaker hypotheses. As in \cite{AmaGoGue} we can
apply Helly's compactness theorem to get a subsequence $u^{h_i}$
converging to some function $u$ in $\L1_{loc}$ whose total variation
in space is uniformly bounded for all $t\ge 0$. Moreover, working as
in \cite[Proposition 5.1]{am2}, one can prove that $u^{h_i}(t,\cdot)$
converges in $\L1$ to $u(t,\cdot)$, for all $t\ge 0$.

\begin{Theorem}\label{th:limit_in_h}
  Let $u^{h_i}$ be a subsequence of solutions of equation
  \dref{eq:hequation} with uniformly bounded total variation
  converging as $i \rightarrow +\infty$ in $\L1$ to some function $u$.
  Then $u$ is a weak solution to the Cauchy problem
  \dref{integrable:eq:bal}--\dref{init:data}.
\end{Theorem}

We omit the proofs of Theorem \ref{first-convergence} and \ref{th:limit_in_h}
since they are very similar to the proofs of \cite[Theorem 3 and 4]{AmaGoGue}. We only observe that, in those proofs, the computations which rely on the $\L\infty$ bound on the source term have to be substituted by the following estimates.
\begin{itemize}
 \item Concerning the proof of Theorem \ref{first-convergence}:
 \begin{displaymath}
  \begin{split}
   &\int_0^h\left|g\left(jh+s,u^\varepsilon(t,jh-)\right)-g\left(jh+s,u(t,jh-)\right)\right|\;ds\\
   &\qquad\le\int_0^h\left\|g(jh+s,\cdot)\right\|_{\C1}\cdot\left|
   u^\varepsilon(t,jh-)-u(t,jh-)\right|\;ds\\
   &\qquad\le\tilde\varepsilon_h\cdot \left|
   u^\varepsilon(t,jh-)-u(t,jh-)\right|.
  \end{split}
 \end{displaymath}
\item Concerning the proof of Theorem \ref{th:limit_in_h}:
 \begin{displaymath}
  \begin{split}
   &\int_0^h\left|g\left(jh+s,u^h(t,jh-)\right)\right|\;ds
   \le \int_0^h\left\|g(jh+s,\cdot)\right\|_{\C1}\;ds\le\tilde\varepsilon_h
  \end{split}
 \end{displaymath}
 and
  \begin{displaymath}
  \begin{split}
   &\int_0^h\left|g\left(jh+s,u^h(t,jh-)\right)-g\left(jh+s,u(t,jh+s)\right)\right|\;ds\\
   &\qquad\le\int_0^h\left\|g(jh+s,\cdot)\right\|_{\C1}\cdot\left|
   u^h(t,jh-)-u(t,jh+s)\right|\;ds\\
   &\qquad\le\tilde\varepsilon_h\cdot \tv\left\{
   u^h(t,\cdot),[(j-1)h,(j+1)h]\right\}\\
   &\qquad\quad+\int_{jh}^{(j+1)h}\omega(x)\left|u^h(t,x)-u(t,x)\right|.
  \end{split}
 \end{displaymath}

\end{itemize}

We observe that all the computations done in \cite[Section 4]{AmaGoGue} rely on the source $g$ only through the amplitude of the zero waves
and on the interaction estimates. Therefore the following two theorems still hold in the more general setting.

\begin{Theorem}\label{tsoin-tsoin}
There exists $\d>0$ such that if $\left\|\om\right\|_{\L1 (\reali)}$
is sufficiently
small, then for any (small) $h>0$ there exist a non empty closed domain
$\cD_h(\d)$ and a unique uniformly
Lipschitz semigroup
$P^h:[0,+\infty)\times\cD_h(\d)\rightarrow\cD_h(\d)$ whose
trajectories $u(t,.)=P_t^h u_o$ solve
\dref{eq:hequation} and are obtained as limit of any sequence of
$\eps,h$--approximate solutions as $\eps$ tends to zero with fixed
$h$.
In particular the semigroups $P^h$ satisfy for any
$u_o,v_o\in\cD_h(\d)$, $t,s\ge 0$
\be
P^h_0 u_o=u_o,\qquad P^h_t \circ P^h_su_o=
P^h_{s+t}u_o,
\label{eq:semi_prop_Ph}
\ee
\be
\left\|P^h_t u_o-P^h_s v_o
\right\|_{\L1(\reali)} \leq
L\left[\left\|u_o-v_o\right\|_{\L1(\reali)}+|t-s|\right]
\label{eq:lip_prop_Ph}
\ee 
for some $L>0$, independent on $h$.
\end{Theorem}

\begin{Theorem}\label{th:final_existence}
If $\left\|\om\right\|_{\L1 (\reali)}$ is sufficiently small, there exist a
constant $L>0$, a non empty closed domain $\cD$ of integrable
functions with small total variation and a semigroup
$P:[0,+\infty)\times \cD\rightarrow\cD$ with the following properties
\begin{enumerate}[i)]
\item
$P_0u=u,\quad\forall u\in\cD;\qquad P_{t+s}u=P_t \circ P_su,\quad\forall
u\in\cD,\ t,s\ge 0.$
\item
$\|P_su-P_tv\|_{\L1 (\reali)}\le L\Big(|s-t|+\|u-v\|_{\L1
(\reali)}\Big),\quad\forall
u\in\cD,\ t,s\ge 0.$
\item
for all $u_o\in\cD$, the function $u(t,\cdot)=P_t u_o$ is a weak
entropy solution of system \dref{integrable:eq:bal}.
\item
for some $\d>0$ and all $h>0$ small enough $\cD\subset\cD_h(\d)$.
\item
There exists a sequence of semigroups $P^{h_i}$ such that $P^{h_i}_t
u$ converges in $\L1$ to $P_tu$ as $i\rightarrow +\infty$ for any
$u\in\cD$.
\end{enumerate}
\end{Theorem}
\begin{Remark}
 \label{newrem}
 Looking at \cite[(4.6)]{AmaGoGue} and the proof of \cite[Theorem 7]{AmaGoGue} one realizes that the invariant domains $\mathcal{D}_h(\delta)$ and $\mathcal{D}$ depend on the particular source term $g(x,u)$. On the other hand estimate \cite[(4.4)]{AmaGoGue} shows that all these domains contain all integrable functions with sufficiently small total variation. Since the bounds $\O(1)$ in Lemma \ref{lemma:equiv_tot_var} depend only on $f$ and on $\|\tilde M\|_{\L1}$, also the constant $C_1$ in \cite[(4.4)]{AmaGoGue} depends only on $f$ and on $\|\tilde M\|_{\L1}$. Therefore there exists $\tilde\delta>0$ depending only on $f$ and on $\|\tilde M\|_{\L1}$ such that $\mathcal{D}_h(\delta)$ and $\mathcal{D}$ contain all integrable functions $u(x)$ with $\tv\left\{u\right\}\le \tilde\delta$. 
\end{Remark}

\section{Uniqueness of $\BV$ entropy solutions}
\label{sec:uniqueness}

The proof of uniqueness in \cite{AmaGoGue} strongly depends on the boundedness of the source, therefore we have to consider it in a more careful way.

\subsection{Some preliminary results}

As in \cite{AmaGoGue} we shall make
use of the following technical lemmas whose proofs can be found in
\cite{Bdue}.

\begin{Lemma}\label{lip_on_intervals}
Let $(a,b)$ a (possibly unbounded) open interval, and let
$\hat\lambda$ be an upper bound for all wave speeds. If $\bar u,\ \bar
v\in\cD_h(\d)$ then for all $t\ge 0$ and
$h>0$, one has
\be
\label{eq:lip_on_intervals}
\int_{a+\hat\lambda t}^{b-\hat\lambda t}
\left|\left(P^h_t\bar u\right)(x)-\left(P^h_t\bar v\right)(x)\right|dx
\le L\int_a^b
\left|\bar u(x)-\bar v(x)\right|dx.
\ee
\end{Lemma}

\begin{Lemma}\label{distance_on_intervals}
Given any interval $I_0=[a,b]$, define the interval of determinacy
\be\label{eq:int_of_det}
I_t=[a+\hat\lambda t,b-\hat\lambda t],\qquad t<
\frac{b-a}{2\hat\lambda}.
\ee
For every Lipschitz continuous map $w:[0,T]\mapsto\cD_h(\d)$ and $h> 0$:
\be\label{eq:distance_on_intervals}
&&\left\|w(t)-P^h_tw(0)\right\|_{\L1(I_t)} \\ \non
&\le&L\int_0^t \left\{\liminf_{\eta\rightarrow 0}
\frac{ \left\|w(s+\eta)-P^h_\eta w(s)
\right\|_{\L1(I_{s+\eta})}}{\eta}\right\}ds.
\ee
\end{Lemma}

\begin{Remark}\label{rem:salvezza}
Lemmas \ref{lip_on_intervals}, \ref{distance_on_intervals}
hold also substituting $P^h$ with the
operator $P$. In this case we have obviously to substitute the domains
$\cD_h(\d)$ with the domain $\cD$ of Theorem \ref{th:final_existence}.
\end{Remark}

Let now $u_\ell,u_r$ be two nearby states and $\lambda<\hat\lambda$;
we consider the function
\be\label{eq:def_di_v}
v(t,x)=\left\{
\begin{array}{lll}
u_\ell&\mbox{ if }&x<\lambda t+x_o\\
u_r&\mbox{ if }&x\ge\lambda t+x_o.
\end{array}
\right.
\ee

\begin{Lemma}\label{lemma_stime_integrali_omogenee}
Call $w(t,x)$ the self-similar solution given by
the standard homogeneous Riemann Solver
with the Riemann data \dref{eq:riemann_data}.
\begin{enumerate}[\it (i)]
\item
In the general case, one has
\be
\label{eq:stima_senza_ipotesi}
\frac 1 t\int_{-\infty}^{+\infty}
\left|v(t,x)-w(t,x)\right|dx=\O(1)|u_\ell-u_r|;
\ee
\item Assuming the additional relations $u_r=R_i(\sigma)(u_\ell)$ and
  $\lambda=\lambda_i(u_r)$ for some $\sigma >0$, $i=1,\ldots,n$ one
  has the sharper estimate \be
\label{eq:stima_rarefazioni}
\frac 1 t\int_{-\infty}^{+\infty}
\left|v(t,x)-w(t,x)\right|dx=\O(1)\sigma^2;
\ee
\item
Let $u^*\in\Omega$ and call $\lambda_1^*<\ldots<\lambda_n^*$ the
eigenvalues of the matrix $A^*=\nabla f(u^*)$. If for some $i$ it holds
$A^*(u_r-u_\ell)=\lambda_i^*(u_r-u_\ell)$ and $\lambda=\lambda_i^*$ in
\dref{eq:def_di_v}
then one has
\be
\label{item:stima_per_lin}
\frac 1 t\int_{-\infty}^{+\infty}
\left|v(t,x)-w(t,x)\right|dx 
=\O(1)|u_\ell-u_r|\Big(
|u_\ell-u^*|+|u^*-u_r|\Big);
\ee
\end{enumerate}
\end{Lemma}

We now prove the next result which is directly related to
our $h$-Riemann solver.

\begin{Lemma}\label{lemma_stime_integrali_sorgente}
Call $w(t,x)$ the self-similar solution given by
the $h$--Riemann Solver in $x_o$
with the Riemann data \dref{eq:riemann_data}.
\begin{enumerate}[\it (i)]
\item
\label{item:hstima_senza_ipotesi}
In the general case one has
\be
\label{eq:hstima_senza_ipotesi}
\frac 1 t\int_{-\infty}^{+\infty}
\left|v(t,x)-w(t,x)\right|dx=\O(1)\Big(|u_\ell-u_r|+\int_0^h\omega(x_o+s)\;ds\Big);
\ee
\item
\label{item:hstima_per_lin}
Assuming the additional relation 
\begin{displaymath}
u_r=u_\ell+\left[\nabla
  f\right]^{-1}(u^*)\displaystyle{\int_{0}^{h}}g(x_o+s,u^*)ds
\end{displaymath}

 with
$\lambda = 0$ in \dref{eq:def_di_v} one has the sharper estimate 
\be\nonumber
&&\frac 1 t\int_{-\infty}^{+\infty}
\left|v(t,x)-w(t,x)\right|dx \\
&=&\O(1)
\Big(\int_0^h\omega(x_o+s)\;ds+|u_\ell - u^*|\Big)
\cdot\int_0^h\omega(x_o+s)\;ds.
\label{eq:hstima_per_lin}
\ee
\end{enumerate}
\end{Lemma}

{\it Proof.}  Estimate {\it \dref{item:hstima_senza_ipotesi}} is a direct
consequence of Lemma \ref{lemma:equiv_tot_var}. 
Let us prove now {\it \dref{item:hstima_per_lin}}. 
Since $\lambda=0$ we derive
\be\non
&&\frac 1 t\int_{-\infty}^{+\infty}
\left|v(t,x)-w(t,x)\right|dx\\
&=&\frac 1 t\int_{-\hat\lambda t}^{0}\left|
u_\ell-w(t,x)\right|dx 
 +\frac 1 t\int_{0}^{\hat\lambda t}\left|
u_r-w(t,x)\right|dx\\
&=&\O(1)\left[\sum_{\iota=1}^p|\sigma_\iota|+\sum_{\iota=p+1}^n|\sigma_\iota|
\right]=\O(1)|\ss|.
\non
\ee
This leads to
\be\non
|\ss|&=&\Big|E\left[h,u_\ell,u_r;x_o\right]-
E\left[h,u_\ell,\Phi_h(x_o,u_\ell);x_o\right]
\Big|\\
&=&\O(1)\left|u_r-\Phi_h(x_o,u_\ell)\right|.
\non
\ee
To estimate this last term, we define $b(y,u)=f^{-1}\left(f(u)+y\right)$
and compute
for some $y_1,\; y_2$:
\be\non
\left|u_\ell+\left[\nabla f\right]^{-1}(u^*)y_1-b(y_2,u_\ell)
\right|&\le& \O(1)|y_1|\cdot |u^*-u_\ell|+\O(1)|y_1-y_2|
\\\non
&&+\left|u_\ell+\left[\nabla f\right]^{-1}(u_\ell)y_2
-b\left(y_2,u_\ell\right)\right|.
\ee
The function $z(y_2)=u_\ell+\left[\nabla f\right]^{-1}(u_\ell)y_2
-b\left(y_2,u_\ell\right)$ satisfies $z(0)=0$, $D_{y_2}z(0)=0$, 
hence we have the estimate
\be\non
\left|u_\ell+\left[\nabla f\right]^{-1}(u^*)y_1-b(y_2,u_\ell)
\right|\le \O(1)\left[|y_1|\cdot |u^*-u_\ell|+
|y_1-y_2|+\left|y_2\right|^2\right].
\ee
If in this last expression we substitute
$$
y_1=\int_{0}^{h}g(x_o+s,u^*)ds,\qquad
y_2=\int_{0}^{h}g(x_o+s,u_\ell)ds
$$
then, we get
$$
\left|u_r-\Phi_h(x_o,u_\ell)\right|=O(1)\Big(\int_0^h\omega(x_o+s)\;ds+|u_\ell-u^*|\Big)\int_0^h\omega(x_o+s)\;ds
$$
which proves \dref{eq:hstima_per_lin}. \square

\subsection{Characterization of the trajectories of $P$}
\label{section4.2}

In this section we are about to give necessary and sufficient conditions for
a function $u(t,\cdot)\in \cD$ to coincide with a semigroup's trajectory.
To this end, we prove the uniqueness of the semigroup $P$ and the
convergence of all the sequence of semigroups $P^h$
towards $P$ as $h \rightarrow 0$.

We begin by introducing some notations: given a $BV$ function $u=u(x)$
and a point $\xi\in\reali$, we denote by $U^\sharp_{(u;\xi)}$ the
solution of the homogeneous Riemann Problem \dref{eq:riemann_data}
with data
\be
u_\ell=\lim_{x\rightarrow\xi -}u(x),\qquad
u_r=\lim_{x\rightarrow\xi +}u(x),\qquad x_o=\xi.
\label{eq:def_riemann_data}
\ee
Moreover we define $U^\flat_{(u;\xi)}$ as the solution of the
linear hyperbolic Cauchy problem with constant coefficients
\be
w_t+\widetilde Aw_x=\widetilde g(x), \qquad w(0,x)=u(x),
\label{eq:def_linear_problem}
\ee
with $\widetilde A = \nabla f \left(u(\xi)\right)$, $\widetilde
g(x)=g\left(x,u(\xi)\right)$.

We will need also the following approximations of $U^\flat_{(u;\xi)}$. 
Let $\overline v$ be a piecewise constant function. We will call $w^h$ the solution of the following Cauchy problem:
\begin{displaymath}
(w^h)_t+\widetilde A(w^h)_x=\sum_{j\in\interi}\d(x-jh)
\int_{0}^{h}\widetilde g(jh+s)\;ds,
\qquad w^h(0,x)=\overline{v}(x).
\end{displaymath}

Define $u^*\dot= u(\xi)$ and let
$\lambda_i=\lambda_i(u^*)$, $r_i=r_i(u^*)$, $l_i=l_i(u^*)$ be respectively
the $i^{th}$ eigenvalue, the $i^{th}$ right/left eigenvectors of
the matrix $\widetilde A$.  As in \cite{AmaGoGue} $w$ and $w^h$ have the following explicit representation 
\begin{eqnarray}\non
w(t,x)&=&\sum_{i=1}^n\left\{\left\langle l_i,u\left(x-\lambda_it\right)\right\rangle +
\frac 1{\lambda_i}\int_{x-\lambda_it}^x
\left\langle l_i,\widetilde g(x')\right\rangle dx'
\right\} r_i\\
\label{eq:happrox_lin_sol}
w^h(t,x)&=&\sum_{i=1}^n\left\{\left\langle l_i,\overline{v}\left(
x-\lambda_it\right)\right\rangle +
\frac 1{\lambda_i}
\left\langle l_i,G^h(t,x)
\right\rangle
\right\} r_i,
\end{eqnarray}
where the function $\displaystyle{G^h(t,x)=\sum_{i=1}^nG_i^h(t,x)r_i}$ is
defined by
\be
\label{eq:definition_Gh}
G_i^h(t,x)=
\left\{
\begin{array}{lll}
\displaystyle{
\sum_{j:\; jh\in\left(x-\lambda_it,x\right)}
\int_{0}^{h}\left\langle l_i,\widetilde g
\left(jh+s\right)\right\rangle ds}&\hbox{ if }\lambda_i>0\\
\displaystyle{
-\sum_{j:\; jh\in\left(x,x-\lambda_it\right)}
\int_{0}^{h}\left\langle l_i,
\widetilde g\left(jh+s\right)\right\rangle ds}&\hbox{ if }\lambda_i<0.
\end{array}
\right.
\ee
Using (\ref{epsilonh}) we can compute
\be
\left|G_i^h(t,x)-\int_{x-\lambda_it}^x
\left\langle l_i,\widetilde g(x')\right\rangle dx\right| = \O(1)\tilde\varepsilon_h.
\ee
Hence, for any $a,b\in\reali$ with $a<b$, we have the error estimate
\begin{equation}
\label{eq:distance_between_w_wh}
\int_{a+\hat\lambda t}^{b-\hat\lambda t}
\left|w(t,x)-w^h(t,x)\right|dx
\le \O(1)\left[\int_{a}^{b}
\left|u(x)-\overline{v}(x)\right|dx +(b-a)\tilde\varepsilon_h\right].
\end{equation}
From \dref{eq:happrox_lin_sol}, \dref{eq:definition_Gh}, it is easy to
see that $w^h(t,x)$ is piecewise constant with discontinuities
occurring along finitely many lines on compact sets in the $(t,x)$
plane for $t\ge 0$. Only finitely many wave front interactions occur
in a compact set, and jumps can be of two types: contact
discontinuities or zero waves. The zero waves are located at the
points $jh$, $j\in \interi$ and satisfy
\be
w^h(t,jh+)-w^h(t,jh-)=\left[ \nabla f\right]^{-1}(u^*)
\int_{jh}^{(j+1)h}\widetilde g\left(jh+s\right)ds.
\label{eq:abra}
\ee
Conversely a contact discontinuity of the $i^{th}$ family located at
the point $x_\al(t)$ satisfies $\dot x_\al(t)=\lambda_i(u^*)$ and
\be
w^h(t,x_\al(t)+)-w^h(t,x_\al(t)-)=\sigma r_i(u^*)
\label{eq:cadabra}
\ee
for some $\sigma\in\reali$.

Now, we can state the uniqeness result in our more general setting.

\begin{Theorem}\label{th:characterisation}
Let $P:\cD\times [0,+\infty)\rightarrow \cD$ be the semigroup
of Theorem
\ref{th:final_existence} and let $\hat\lambda$ be an upper bound  for
all wave speeds. Then every trajectory $u(t,\cdot)=P_t u_{0}$,
$u_{0} \in \cD$, satisfies
the following conditions at every $\tau \ge 0$.
\begin{enumerate}[(\it i)]
\item
\label{first_int_cond}
For every $\xi$, one has
\be
\label{eq:first_int_cond}
\lim_{\theta\rightarrow 0}
\frac 1 {\theta}\int_{\xi-\theta\hat\lambda}^{\xi+\theta\hat\lambda}
\left|u(\tau+\theta,x)-U^\sharp_{(u(\tau);\xi)}\left(\theta,x\right)\right|dx
=0.
\ee
\item
\label{second_int_cond}
There exists a constant $C$ such that, for every $a<\xi<b$ and
$0<\theta<\frac{b-a}{2\hat\lambda}$, one has
\be
\label{eq:second_int_cond}
\frac 1 {\theta} \int_{a+\theta\hat\lambda}^{b-\theta\hat\lambda}
\left|u(\tau+\theta,x)-U^\flat_{(u(\tau);\xi)}
\left(\theta,x\right)\right| dx \\ \non
\le  C\Big[\tv\left\{u(\tau);(a,b)\right\} +\int_a^b\omega(x)\;dx\Big]^2.
\ee
\end{enumerate}

Viceversa let $u:[0,T]\rightarrow\cD$ be Lipschitz continuous as a map
with values in $\L1(\reali,\reali^n)$ and assume that the conditions
\dref{first_int_cond}, \dref{second_int_cond} hold at almost every
time $\tau$. Then
$u(t,\cdot)$ coincides with a trajectory of the semigroup $P$.
\end{Theorem}

\begin{Remark}
\label{rem:recall}
 The difference with respect to the result in \cite{AmaGoGue} is the presence of the integral in the right hand side of formula (\ref{eq:second_int_cond}). If $\omega$ is in $\L\infty$, the integral can be bounded by $\O(1)(b-a)$ and we recover the estimates in \cite{AmaGoGue}. Note also that the quantity 
\begin{displaymath}
 \mu\left((a,b)\right)=\tv\left\{u(\tau);(a,b)\right\} +\int_a^b\omega(x)\;dx
\end{displaymath}
is a uniformly bounded finite measure and this is what is needed for proving the sufficiency part of the above Theorem.
\end{Remark}

\begin{proof}
\textbf{Part 1: Necessity} Given a semigroup trajectory
$u(t,\cdot)=P_t\bar u$, $\bar u\in\cD$ we now show that the conditions
{\it \dref{first_int_cond}},
{\it \dref{second_int_cond}} hold for every $\tau\ge 0$.

As in \cite{AmaGoGue} we use the following notations. For fixed 
$h,\;\theta,\;\ve>0$ we define
$J_t=J_t^-\cup J^o_t\cup J^+_t$
with
\be\non
J_t^-&=&\left(\xi-\left(2\theta-t+\tau\right)\hat \lambda,\xi-(t-\tau
)\hat\lambda\right);\\\label{intervals}
J_t^o&=&\left[\xi-(t-\tau)\hat\lambda,\xi+(t-\tau
)\hat\lambda\right];\\\non
J_t^+&=&\left(\xi+(t-\tau
)\hat\lambda,\xi+\left(2\theta-t+\tau\right)\hat \lambda\right).
\ee
Let $U^{\sharp,\ve}_{(u(\tau);\xi)}\left(\theta,x\right)$
be the piecewise constant function
obtained from $U^{\sharp}_{(u(\tau);\xi)}\left(\theta,x\right)$
dividing the centered rarefaction waves in equal parts and replacing them by
rarefaction fans containing wave fronts whose strength is less than $\ve$.
Observe that:
\be
\frac 1
t\int_{-\infty}^{+\infty}\left|U^{\sharp,\ve}_{(u(\tau);\xi)}\left(\theta,x\right
)
-U^{\sharp}_{(u(\tau);\xi)}\left(\theta,x\right)\right|\;dx=\O(1)\ve.
\label{eq:compare_to_eps}
\ee
Applying estimate \dref{eq:distance_on_intervals} to the function $U^{\sharp,\ve}_{(u(\tau);\xi)}$
we obtain
\be
\label{eq:stimapiu}
&&\int_{J_{\tau+\theta}}\left|
U^{\sharp,\ve}_{(u(\tau);\xi)}(\theta,x)-\left(P^h_\theta U^{\sharp,\ve}_{(u(\tau);\xi)}(0)\right)(x)
\right|dx \\ \non
&\le& L
\int_\tau^{\tau+\theta}\liminf_{\eta\rightarrow 0}
\frac{\left\|U^{\sharp,\ve}_{(u(\tau);\xi)}(t-\tau+\eta)-P_\eta^hU^{\sharp,\ve}_{(u(\tau);\xi)}(t-\tau)
\right\|_{\L1\left(J_{t+\eta}\right)}}
{\eta}dt.
\ee
The discontinuities of $U^{\sharp,\ve}_{(u(\tau);\xi)}$ do not cross the Dirac comb for almost
all times $t\in (\tau,\tau+\theta)$. Therefore we compute for such a time
$t$:
\be
\label{eq:uno_asterisco}
&&\frac 1{\eta} \int_{J_{t+\eta}}
\left|U^{\sharp,\ve}_{(u(\tau);\xi)}(t-\tau+\eta,x)-\left(P_\eta^h U^{\sharp,\ve}_{(u(\tau);\xi)}(t-\tau)
\right)(x)\right|dx\\\non
&=&\frac 1{\eta}
\int_{J^-_{t+\eta}\cup J^o_{t+\eta}\cup J^+_{t+\eta}}
\left|U^{\sharp,\ve}_{(u(\tau);\xi)}(t-\tau+\eta,x)-\left(P_\eta^hU^{\sharp,\ve}_{(u(\tau);\xi)} (t-\tau)
\right)(x)\right|dx.
\ee
Define ${\cal{W}}_t$ the set of points in which $U^{\sharp,\ve}_{(u(\tau);\xi)}(t-\tau)$ has a discontinuity while ${\cal{Z}}_h$ is the set of points in which the zero waves are located. 
If $\eta$ is sufficiently small, the solutions of the Riemann problems arising at the discontinuities of $U^{\sharp,\ve}_{(u(\tau);\xi)}(t-\tau)$ do not interact, therefore
\begin{eqnarray}\non
  &   &\frac 1{\eta} \int_{J^o_{t+\eta}}
\left|U^{\sharp,\ve}_{(u(\tau);\xi)}(t-\tau+\eta,x)-\left(P_\eta^hU^{\sharp,\ve}_{(u(\tau);\xi)}(t-\tau) 
\right)(x)\right|dx \\\non
  & = &\left(\sum_{x\in J^o_t\cap{\cal{W}}_t}+\sum_{x\in J^o_t\cap{\cal{Z}}_h}\right)\\\non
&& \frac 1{\eta} \int_{x-\hat\lambda\eta}^{x+\hat\lambda\eta}
\left|U^{\sharp,\ve}_{(u(\tau);\xi)}(t-\tau+\eta,y)-\left(P_\eta^hU^{\sharp,\ve}_{(u(\tau);\xi)}(t-\tau)
\right)(y)\right|dy
\end{eqnarray}
Note that the shock are solved exactly both in $U^{\sharp,\ve}_{(u(\tau);\xi)}$ and in $P^hU^{\sharp,\ve}_{(u(\tau);\xi)}$ therefore they make no contribution in the summation. 
To estimate the approximate rarefactions we use the estimate (\ref{eq:stima_rarefazioni}) hence
\begin{eqnarray}\non
&  & \sum_{x\in J^o_t\cap{\cal{W}}_t}
\frac 1{\eta} \int_{x-\hat\lambda\eta}^{x+\hat\lambda\eta}
\left|U^{\sharp,\ve}_{(u(\tau);\xi)}(t-\tau+\eta,x)-\left(P_\eta^hU^{\sharp,\ve}_{(u(\tau);\xi)} (t-\tau)
\right)(x)\right|dx\\
&&\le 
\O(1)\sum_{\stackrel{x\in J^o_t\cap{\cal{W}}_t}{rarefaction}}|\sigma|^2\label{eq:49}
 \le \O(1)\varepsilon \tv\left\{U^{\sharp,\ve}_{(u(\tau);\xi)}(t-\tau);J^0_{t}\right\}\\\non
&&\le\O(1)\varepsilon\left|u(\tau,\xi+)-u(\tau,\xi-)\right|
\end{eqnarray}
Concerning the zero waves, recall that $t$ is chosen such that $U^{\sharp,\ve}_{(u(\tau);\xi)}$ is constant there, and $P^h$ is the exact solution of an $h$--Riemann problem, hence we can apply (\ref{eq:hstima_senza_ipotesi}) with $u_\ell=u_r$ and
obtain
\begin{eqnarray}\non
&  & \sum_{x\in J^o_t\cap{\cal{Z}}_h}
\frac 1{\eta} \int_{x-\hat\lambda\eta}^{x+\hat\lambda\eta}
\left|U^{\sharp,\ve}_{(u(\tau);\xi)}(t-\tau+\eta,x)-\left(P_\eta^hU^{\sharp,\ve}_{(u(\tau);\xi)} (t-\tau)
\right)(x)\right|dx\\\label{eq:51}
&\le &
\O(1)\sum_{jh\in J^o_t}\int_{0}^h\omega(jh+s)\;ds
\le\O(1)\left(\int_{J^o_t}\omega(x)\;dx + \tilde\varepsilon_h\right)
\end{eqnarray}
Finally using (\ref{eq:51}) and (\ref{eq:49}) we get in the end
\be
\label{eq:due_asterisco}
\frac 1{\eta} \int_{J^o_{t+\eta}}
\left|U^{\sharp,\ve}_{(u(\tau);\xi)}(t-\tau+\eta,x)-\left(P_\eta^hU^{\sharp,\ve}_{(u(\tau);\xi)} (t-\tau)
\right)(x)\right|dx \\ \non =\O(1)\left\{\int_{J^0_t}\omega(x)\;dx+\tilde\varepsilon_h+\ve\right\}.
\ee
Moreover, following the same steps as before and using
\dref{eq:stima_senza_ipotesi} and \dref{eq:hstima_senza_ipotesi} with
$u_\ell=u_r$ we get
\be
\label{eq:tre_asterisco}
\frac 1{\eta} \int_{J^+_{t+\eta}}
\left|U^{\sharp,\ve}_{(u(\tau);\xi)}(t-\tau+\eta,x)-\left(P_\eta^hU^{\sharp,\ve}_{(u(\tau);\xi)} (t-\tau)
\right)(x)\right|dx \\=  \non 
\O(1) \left\{\int_{J^+_t}\omega(x)\;dx+\tilde\varepsilon_h
\right\}.
\ee
Note that here there is no total variation of $U^{\sharp,\ve}_{(u(\tau);\xi)}$ since in $J^+_t$ it is constant. 
A similar estimate holds for the interval $J^-_{t+\eta}$.
Putting together
\dref{eq:uno_asterisco}, \dref{eq:due_asterisco}, \dref{eq:tre_asterisco},
one has
\be\non
\frac 1{\eta} \int_{J_{t+\eta}}
\left|U^{\sharp,\ve}_{(u(\tau);\xi)}(t-\tau+\eta,x)-\left(P_\eta^h U^{\sharp,\ve}_{(u(\tau);\xi)}(t-\tau)
\right)(x)\right|dx \\=  \non \O(1)\Big(
\int_{J_\tau}\omega(x)\;dx+\tilde\varepsilon_h+\ve
\Big).
\ee
Hence, setting $\tilde v=U^{\sharp,\ve}_{(u(\tau);\xi)}(0)=U^{\sharp}_{(u(\tau);\xi)}(0)$ by \dref{eq:stimapiu}, 
we have
\be
\label{eq:stimameno}
\int_{J_{\tau+\theta}}\left|
U^{\sharp,\ve}_{(u(\tau);\xi)}(\theta,x)-\left(P^h_\theta \tilde v\right)(x)
\right|dx=\O(1)\theta\Big(
\int_{J_\tau}\omega(x)\;dx+\tilde\varepsilon_h+\ve
\Big).
\ee
Finally we take the sequence $P^{h_i}$ converging to $P$.
Using \dref{eq:lip_on_intervals} we have
\be\label{after51}
\frac{1}{\theta}\left\|P^{h_i}_\theta u(\tau)-P^{h_i}_\theta \tilde v
\right\|_{\L1\left(J_{\tau+\theta}\right)}&\le& \frac{1}{\theta}L\left\|u(\tau)-\tilde v
\right\|_{\L1\left(J_{\tau}\right)}\\\non
&=&\frac{L}{\theta}\int_{\xi-2\hat\lambda\theta}^{\xi+2\hat\lambda\theta}\left|u(\tau,x)-\tilde v(x)\right|dx\\ \non &\dot=&\bar\varepsilon_\theta,
\ee
where $\bar\varepsilon_{\theta}$ tends to zero as $\theta$ tends to zero due to the fact that $u(\tau)$ has right and left limit at any point: for any given $\epsilon>0$ if $\theta$ is sufficiently small 
$|u(\tau,x)-\tilde v(x)|=|u(\tau,x)-u(\tau,\xi-)|\le\epsilon$ for $x\in (\xi-2\hat\lambda\theta,\xi)$.

Therefore by \dref{eq:compare_to_eps}, \dref{eq:stimameno}, we
derive:
\be\non
&&\frac 1 {\theta}\int_{\xi-\theta\hat\lambda}^{\xi+\theta\hat\lambda}
\left|u(\tau+\theta,x)-U^\sharp_{(u(\tau);\xi)}\left(\theta,x\right)\right|dx
\\ &=& \non \frac{\left\|
P_\theta u(\tau)-P^{h_i}_\theta u(\tau)
\right\|_{\L1(\reali)}}{\theta}+\bar\varepsilon_\theta
+\O(1)\left[
\int_{J_\tau}\omega(x)\;dx+\tilde\varepsilon_{h_i}
\right].
\ee
The left hand side of the previous estimate does not depend on $\ve$
and $h_i$, hence
\be\non
\frac 1 {\theta}\int_{\xi-\theta\hat\lambda}^{\xi+\theta\hat\lambda}
\left|u(\tau+\theta,x)-U^\sharp_{(u(\tau);\xi)}\left(\theta,x\right)\right|dx
=  O(1)
\int_{J_\tau}\omega(x)\;dx+\bar\varepsilon_{\theta}.
\ee
Note that the intervals $J_\tau$ depend on $\theta$ (see \ref{intervals}). So taking the limit as $\theta\rightarrow 0$ in the previous estimate
yields \dref{eq:first_int_cond}.
\vspace{0.5truecm}

To prove {\it \dref{second_int_cond}} let $\theta>0$ and a point
$(\tau,\xi)$ be given together with an open interval $(a,b)$
containing $\xi$. Fix $\ve>0$ and choose a piecewise constant function $\bar
v\in\cD$ satisfying $\bar v(\xi)=u(\tau,\xi)$ together with
\begin{equation}
\int_a^b\left|\bar v(x)-u(\tau,x)\right|\;dx\le\ve,
\quad
\tv\left\{
\bar v;(a,b)
\right\}
\le
\tv\left\{
u(\tau);(a,b)
\right\}
\label{eq:second_totvar_prop}
\end{equation}
Let now $w^h$ be defined by \dref{eq:happrox_lin_sol} ($u^*=\bar v(\xi)=u(\tau,\xi)$). From
\dref{eq:distance_between_w_wh},
\dref{eq:second_totvar_prop} we have the estimate
\be
\label{eq:propp}
\int_{a+\theta\hat\lambda}^{b-\theta\hat\lambda}
\left|U^\flat_{(u(\tau);\xi)}\left(\theta,x\right)-
w^h(\theta,x)
\right|dx\le\O(1)\Big(\ve +\tilde\varepsilon_h(b-a)\Big).
\ee
Using \dref{eq:int_of_det}, \dref{eq:distance_on_intervals} we get
\be
\label{eq:prip}
&&\int_{a+\theta\hat\lambda}^{b-\theta\hat\lambda}
\left|w^h(\theta,x)-\left(P_\theta^h w^h(0)\right)(x)
\right|dx\\
&\le & \non L
\int_{\tau}^{\tau +\theta}\liminf_{\eta\rightarrow 0}
\frac{ \left\|w^h(t-\tau+\eta)-P^h_\eta w^h(t-\tau)
\right\|_{\L1(\tilde I_{t+\eta})}}{\eta}dt
\ee
where we have defined $\tilde I_{t+\eta}=I_{t-\tau+\eta}$.
Let $t\in (\tau,\tau+ \theta)$ be a time for which there is no
interaction in $w^h$; in particular,
discontinuities which travel with a non-zero velocity
do not cross the Dirac comb (this happens for almost all $t$).
We observe that by the explicit formula \dref{eq:happrox_lin_sol}:
\begin{equation}
\label{eq:3;35}
\tv\left\{w^h(t-\tau); \tilde I_{t} \right\}=\O(1)
\Big(\tv\left\{\bar v; (a,b) \right\}+\int_a^b\omega(x)\;dx+\tilde\varepsilon_h\Big)
\end{equation}
\begin{equation}
\label{eq:3;36}
\left|w^h(t-\tau,x)-\bar v(\xi)\right|=\O(1)\Big(
\tv\left\{\bar v; (a,b) \right\}+\int_a^b\omega(x)\;dx+\tilde\varepsilon_h\Big).
\end{equation}
As before for $\eta$ sufficiently small we can split homogeneous and zero waves
\begin{eqnarray}
  &   &\frac 1 {\eta}\int_{\tilde I_{t+\eta}}\left|w^h(t-\tau+\eta,x)-\left(
P^h_\eta w^h(t-\tau)
\right)(x)
\right|dx
 \\\non
  & = &\left(\sum_{x\in \tilde I_t\cap{\cal{W}}_t}+\sum_{x\in \tilde I_t\cap{\cal{Z}}_h}\right)
\frac 1{\eta} \int_{x-\hat\lambda\eta}^{x+\hat\lambda\eta}
\left|w^h(t-\tau+\eta,x)-\left(
P^h_\eta w^h(t-\tau)
\right)(x)
\right|dx
\end{eqnarray}

The homogeneous waves in $w^h$ satisfy (\ref{eq:cadabra}), with $\bar v(\xi)$ in place of $u^*$, hence we can apply 
(\ref{item:stima_per_lin}) which together with  \dref{eq:3;35}, (\ref{eq:3;36}) leads to 

\begin{eqnarray}\non
  &  &\sum_{x\in \tilde I_t\cap{\cal{W}}_t}
\frac 1{\eta} \int_{x-\hat\lambda\eta}^{x+\hat\lambda\eta}
\left|w^h(t-\tau+\eta,x)-\left(
P^h_\eta w^h(t-\tau)
\right)(x)
\right|dx\\\non
&\le& \O(1)\sum_{x\in \tilde I_t\cap{\cal{W}}_t}|\Delta w^h(t-\tau,x)|\Big(
\tv\left\{\bar v; (a,b) \right\}+\int_a^b\omega(x)\;dx+\tilde\varepsilon_h\Big)\\\non
&\le& \O(1)\tv\left\{w^h(t-\tau),\tilde I_t\right\}\Big(
\tv\left\{\bar v; (a,b) \right\}+\int_a^b\omega(x)\;dx+\tilde\varepsilon_h\Big)\\
\non
&\le& \O(1)\Big(
\tv\left\{\bar v; (a,b) \right\}+\int_a^b\omega(x)\;dx+\tilde\varepsilon_h\Big)^2
\end{eqnarray}
where $\Delta w^h(t-\tau,x)$ denotes the jump of $w^h(t-\tau)$ at $x$.

The zero waves in $w^h$ satisfy (\ref{eq:abra}), hence we can apply 
(\ref{eq:hstima_per_lin}) which together with (\ref{eq:3;36}) leads to 

\begin{eqnarray}\non
  &  &\sum_{x\in \tilde I_t\cap{\cal{Z}}_h}
\frac 1{\eta} \int_{x-\hat\lambda\eta}^{x+\hat\lambda\eta}
\left|w^h(t-\tau+\eta,x)-\left(
P^h_\eta w^h(t-\tau)
\right)(x)
\right|dx\\\non
&\le& \O(1)\sum_{x\in \tilde I_t\cap{\cal{Z}}_h}\int_0^h\omega(x+s)\;ds\cdot
\Big(
\tv\left\{\bar v; (a,b) \right\}+\int_a^b\omega(x)\;dx+\tilde\varepsilon_h\Big)\\\non
&\le& \O(1)\left(\int_{\tilde I_t}\omega(x)\;dx+\tilde\varepsilon_h\right)\Big(
\tv\left\{\bar v; (a,b) \right\}+\int_a^b\omega(x)\;dx+\tilde\varepsilon_h\Big)\\
\non
&\le& \O(1)\Big(
\tv\left\{\bar v; (a,b) \right\}+\int_a^b\omega(x)\;dx+\tilde\varepsilon_h\Big)^2
\end{eqnarray}

Let now $P^{h_i}$ be the subsequence converging to $P$.
Since $w^h(0)=\bar v$
using
\dref{eq:propp}, \dref{eq:prip},
\dref{eq:second_totvar_prop}, and the last estimates
we get
\be\non
&&\frac 1 {\theta} \int_{a+\theta\hat\lambda}^{b-\theta\hat\lambda}
\left|u(\tau+\theta,x)-U^\flat_{(u(\tau);\xi)}\left(\theta,x\right)\right|
dx\\
&\le&\non 
\frac{\|P_\theta u(\tau)-P_\theta^{h_i}u(\tau)\|_{\L1(\reali)}}{\theta}+L\frac{\|u(\tau)-\bar v\|_{\L1(\reali)}}{\theta}
\\\non
&&+\O(1) \bigg\{\frac{\varepsilon+\tilde\varepsilon_{h_i}\cdot(b-a)}{\theta}+ \Big(
\tv\left\{\bar v; (a,b) \right\}+\int_a^b\omega(x)\;dx+\tilde\varepsilon_{h_i}\Big)^2\bigg\}.
\ee
So for $\ve,\; h_i\rightarrow 0$ we obtain the desired inequality.
\\
\textbf{Part 2: Sufficiency} By Remark \ref{rem:salvezza} we can
apply \dref{eq:distance_on_intervals} to $P$ and hence the proof for
the homogeneous case presented in \cite{Bdue}, which relies on the property recalled in Remark \ref{rem:recall}, can be followed exactly
for our case, hence it will be not repeated here. $\square$

\end{proof}

\paragraph{Proof of Theorem \ref{th:introduction}} 
It is now a direct
consequence of Theorems \ref{th:final_existence} and 
\ref{th:characterisation}.  \square

\section{Proofs related to Section \ref{application}}
\label{section5}

Consider the equation
\begin{displaymath}
 u_t+f(u)_x=a'g(u)
\end{displaymath}
for some $a\in \BV$. Equation (\ref{dynamics}) is comprised in this setting with the substitution $a\mapsto \ln a$. For this kind of equations we consider the exact stationary solutions instead of approximated ones as in (\ref{def:Phi}). Therfore call $\Phi(a,\bar u)$ the solution of the following Cauchy problem:
\begin{equation}\label{ODE}
 \left\{
   \begin{array}{l}
    \frac{d}{da}u(a)=\left[D_uf(u(a))\right]^{-1}g(u(a))\\
    u(0)=\bar u    
   \end{array}
\right.
\end{equation}
If $a$ is sufficiently small, the map $u\mapsto \Phi(a,u)$ satisfies Lemma \ref{lemma:Riemann_prob0}. We call $a$-Riemann problem the Cauchy problem
\begin{equation}
\label{eq:62b}
 \left\{
   \begin{array}{l}
    u_t+f(u)_x=a'g(u)\\
    (a,u)(0,x)=\left\{
    \begin{array}{ll}
     (a^-,u_l)&\mbox{ if }x<0\\
     (a^+,u_r)&\mbox{ if }x>0
     \end{array}
    \right.    
   \end{array}
\right.
\end{equation}
its solution will be the function described in Definition \ref{def:h_R_solver}
using the map $\Phi(a^+-a^-,u^-)$ instead of the $\Phi_h$ in there. Observe that if $a^+=a^-$ the $a$-Riemann solver coincides with the usual homogeneous Riemann solver.

\begin{Definition}
Given a function $u\in\BV$ and two states $a^-$, $a^+$, we define $\bar U^\sharp_{u}\left(t,x\right)$ as the solution of the $a$-Riemann solver (\ref{eq:62b}) with $u_l=u(0-)$ and $u_r=u(0+)$.
\end{Definition}

\paragraph{Proof of Theorem \ref{Theorem2}:} 
Since $\|a_l^\prime\|_{\L1}=|a^+-a^-|$, hypothesis $(P_2)$ is satisfied uniformly with respect to $l$, moreover the smallness of  $|a^+-a^-|$ ensures that the $\L1$ norm of $\omega$ in $(P_3)$ is  small. Therefore the hypotheses of Theorem \ref{th:introduction} are satisfied uniformly with respect to $l$.

Let $P^l$ be the semigroup related with the smooth section $a_l$. By Remark \ref{newrem}, if $\tv\left\{u\right\}$ is sufficiently small, $u$ belongs to the domain of $P^l$ for every $l>0$. Since the total variation of $P^l_tu$ is uniformly bounded for a fixed initial data $u$,  Helly's theorem guarantees that there is a converging subsequence $P^{l_i}_tu$. By a diagonal argument one can show that there is a converging subsequence of semigroups converging to a limit semigroup $P$ defined on an invariant domain (see \cite[Proof of Theorem 7]{AmaGoGue}). 

For the uniqueness we are left to prove the integral estimate
(\ref{eq:first_int_cond}) in the origin with $U^\sharp$ subsituted by $\bar U^\sharp$.

Therefore we have to show that the quantity 
\be
\label{eq:first_int_condb}
\frac 1 {\theta}\int_{-\theta\hat\lambda}^{+\theta\hat\lambda}
\left|u(\tau+\theta,x)-\bar U^\sharp_{u(\tau)}\left(\theta,x\right)\right|dx
\ee
converges to zero as $\theta$ tends to zero. We will estimate (\ref{eq:first_int_condb}) in several steps. First define $\bar v=\bar U^\sharp_{u(\tau)}\left(0,x\right)$ and compute
\begin{equation}
 \frac 1 {\theta}\int_{-\theta\hat\lambda}^{+\theta\hat\lambda}
\left|(P_\theta u(\tau))(x)-(P_\theta \bar v)(x)\left(\theta,x\right)\right|dx\le \bar\epsilon_{\theta}.
\end{equation}
as in (\ref{after51}).
Then we consider the approximating sequence $P^{l_i}$ corresponding to the source term $a_{l_i}$ and the semigroups $P^{l_i,h}$ which converge to $P^{l_i}$ in the sense of Theorem 
\ref{th:final_existence}. Hence we have
\begin{displaymath}
\lim_{i\to\infty}\lim_{h\to 0} \frac 1 {\theta}\int_{-\theta\hat\lambda}^{+\theta\hat\lambda}
\left|(P^{l_i,h}_\theta \bar v)(x)-(P_\theta \bar v)(x)\right|dx=0
\end{displaymath}

\begin{figure}[htpb]
  \centering
  \begin{psfrags}
    \psfrag{u}{$\bar U^{\sharp,\varepsilon,l,h}$} 
    \psfrag{x}{$x$} 
    \psfrag{h}{$h$} 
    \psfrag{t}{$t$}
    \psfrag{ul}{$u(\tau,0-)=u_l$}
    \psfrag{ur}{$u(\tau,0+)=u_r$}
    \psfrag{ut}{$\widetilde U$}
    \psfrag{-l}{$-\frac l2$} 
    \psfrag{l}{$\frac l2$} 
    \psfrag{lambda}{$x=\hat\lambda t$}
    \psfrag{lambdat}{$x=\frac l2+\hat\lambda t$} 
    \includegraphics[width=12cm]{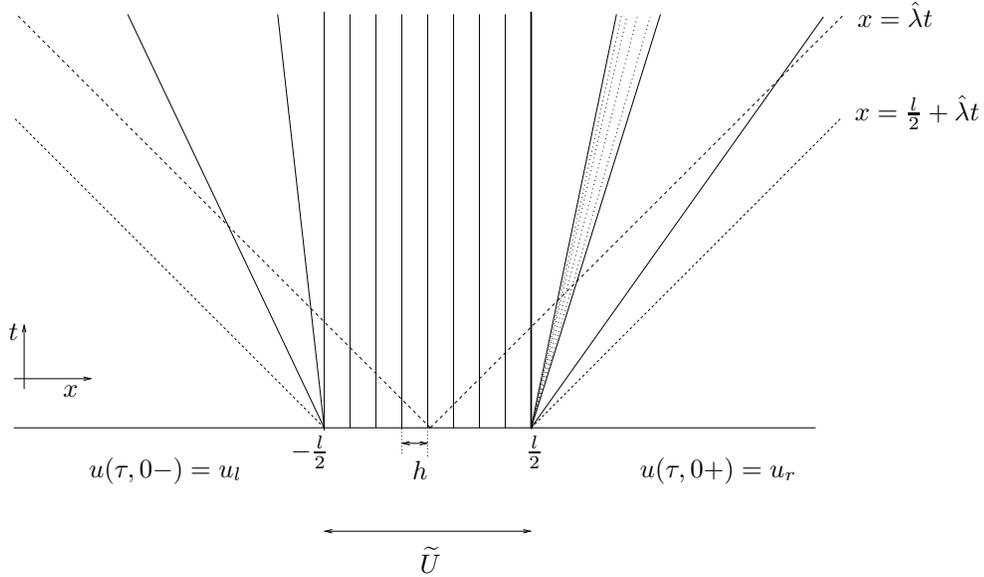}
  \end{psfrags}
  \caption{Illustration of $\bar U^{\sharp,\varepsilon,l,h}$ in the $(t,x)$ plane}
  \label{fig:p}
\end{figure}

For notational convenience we skip the subscript $i$ in $l_i$.
As in (\ref{eq:compare_to_eps}) we approximate rarefactions in $\bar U^\sharp_{u(\tau)}$ introducing the function $\bar U^{\sharp,\varepsilon}_{u(\tau)}$. Then we define (see Figure \ref{fig:p})
\begin{displaymath}
 \bar U^{\sharp,\varepsilon,l,h}_{u(\tau)}(t-\tau,x)=
\left\{
\begin{array}{ll}
 \bar U^{\sharp,\varepsilon}_{u(\tau)}(t-\tau,x+\frac{l}{2})&\text{ for }x< -l/2\\
  \widetilde U(x)&\text{ for }-l/2\le x\le l/2\\
 \bar U^{\sharp,\varepsilon}_{u(\tau)}(t-\tau,x-\frac{l}{2})&\text{ for }x> l/2
\end{array}
\right.
\end{displaymath}
where $\widetilde U(x)$ is piecewise constant with jumps in the points $jh$ satisfying $\widetilde U(jh+)=\Phi(jh,\widetilde U(jh-))$. Furthermore $\widetilde U(-l/2-)= \bar U^{\sharp,\varepsilon}_{(u;\tau)}(t-\tau,0-)$ and $\Phi$ is defined as in (\ref{def:Phi}) using the source term $g(x,u)=a'_l(x)g(u)$. Observe that the jump between $\widetilde U(l/2-)$ and $\bar U^{\sharp,\varepsilon,l,h}_{u(\tau)}(t-\tau,l/2+)$ does not satisfy any jump condition, but as $\widetilde U(x)$ is an ``Euler'' approximation of the ordinary differential equation $f(u)_x=a_{l}^\prime g(u)$, this jump is of order $\tilde\varepsilon_h$.
Since $\bar U^{\sharp,\varepsilon}_{u(\tau)}$ and $\bar U^{\sharp,\varepsilon,l, h}_{u(\tau)}$ have uniformly bounded total variation we have the estimate
\begin{displaymath}
\frac 1 {\theta}\int_{-\theta\hat\lambda}^{+\theta\hat\lambda}
\left|\bar U^{\sharp,\varepsilon}_{u(\tau)}\left(\theta,x\right)-\bar U^{\sharp,\varepsilon,l,h}_{u(\tau)}\left(\theta,x\right)\right|dx\le \O(1)\frac{l}{\theta}
\end{displaymath}
the bound $\O(1)$ not depending on $h$. We apply Lemma \ref{distance_on_intervals} on the remaining term
\begin{eqnarray}\non
&&\frac 1 {\theta}\int_{-\theta\hat\lambda}^{+\theta\hat\lambda}
\left|(P^{l,h}_\theta\bar v)(x)-\bar U^{\sharp,\varepsilon,l,h}_{u(\tau)}\left(\theta,x\right)\right|dx\\\non
&&\le L\int_{\tau}^{\tau+\theta}\liminf_{\eta\to 0}\frac{\|\bar U^{\sharp,\varepsilon,l,h}_{u(\tau)}\left(t-\tau+\eta\right)-P^{l,h}_\eta\bar U^{\sharp,\varepsilon,l,h}_{u(\tau)}\left(t-\tau\right)\|_{\L1(J_{t+\eta})}}{\eta}
\end{eqnarray}
To estimate this last term we proceed as before. Observe that $P^{l,h}$ does not have zero waves outside the interval $[-\frac l2 -h,\frac l2+h]$ since outside the interval $[-\frac l2,\frac l2]$  
 the function $a_l^\prime$ is identically zero. If $\eta$ is small enough, the waves in 
$P^{l,h}_\eta\bar U^{\sharp,\varepsilon,l,h}_{u(\tau)}\left(t-\tau\right)$ do not interact, therefore the 
computation of the $\L1$ norm in the previous integral, as before can be splitted in a summation on 
the points in which there are zero waves in $P^{l,h}$ or jumps in $\bar 
U^{\sharp,\varepsilon,l,h}_{u(\tau)}\left(t-\tau\right)$. Observe that the jumps 
of $\bar U^{\sharp,\varepsilon,l,h}_{u(\tau)}\left(t-\tau+\eta\right)$ in the interval $(-\frac 
l2,+\frac l2)$, are defined exactly as the zero waves in $P^{l,h}$ so
we have no contribution to the summation from this interval. Outside the interval $[-\frac l2 -h,\frac l2+h]$, $P^h$ coincides with the homogeneous semigroup, hence we have only the second order contribution from the approximate rarefactions in $\bar 
U^{\sharp,\varepsilon,l,h}_{u(\tau)}\left(t-\tau\right)$ as in (\ref{eq:49}). Furthermore we might have a zero wave in the interval $[-\frac l2-h,-\frac l2]$ 
and a discontinuity of $\bar U^{\sharp,\varepsilon,l,h}_{u(\tau)}$ in the point $x=\frac l2$ of order $\tilde\varepsilon_h$. Using 
(\ref{eq:hstima_senza_ipotesi}) for the zero wave and (\ref{eq:stima_senza_ipotesi}) for the discontinuity (since $P^h$ is equal to the homogeneous semigroup in $x=\frac l2$), we get
\begin{eqnarray}\non
 \liminf_{\eta\to 0}\frac{\|\bar U^{\sharp,\varepsilon,l,h}_{u(\tau)}\left(t-\tau+\eta\right)-P^{l,h}_\eta\bar U^{\sharp,\varepsilon,l,h}_{u(\tau)}\left(t-\tau\right)\|_{\L1(J_{t+\eta})}}{\eta}\le
\O(1)\left(\varepsilon+\tilde\varepsilon_h\right)
\end{eqnarray}
Which completes the proof if we let first $\varepsilon$ tend to zero, then $h$ tend to zero, then $l$ tend to zero and finally $\theta$ tend to zero.
As in the previous proof, the sufficiency part can be obtained following the proof for
the homogeneous case presented in \cite{Bdue}. 
\square

\paragraph{Proof of Proposition \ref{prop1}:} 
Call $S$ the semigroup defined in \cite{ColomboHertySachers}.
The estimates for this semigroup outside the origin are equal to the ones for the Standard Riemann Semigroup see \cite{Bdue}. Concerning the origin we first observe that the choice (\ref{eq:psieq}) implies that the solution to the Riemann problem in \cite[Proposition 2.2]{ColomboHertySachers} coincides with 
$\bar U^{\sharp}_{u(\tau)}$. We need to show that 
\be
\label{eq:first_int_condbappo}
\lim_{\theta\to0}\frac 1 {\theta}\int_{-\theta\hat\lambda}^{+\theta\hat\lambda}
\left|u(\tau+\theta,x)-\bar U^\sharp_{u(\tau)}\left(\theta,x\right)\right|dx=0.
\ee
with $u(t,x)=(S_tu_o)(x)$.
As before, we first approximate $\bar U^{\sharp}_{u(\tau)}$ with $\bar U^{\sharp,\varepsilon}_{u(\tau)}$ and $u(\tau)$ with $\bar U^{\sharp}_{u(\tau)}(0)\dot =\bar v$ then we apply Lemma \ref{distance_on_intervals} (which holds also for the semigroup $S$) and compute
\be\non
&&\frac 1 {\theta}\int_{-\theta\hat\lambda}^{+\theta\hat\lambda}
\left|(S_\theta \bar v)(x)-\bar U^{\sharp,\varepsilon}_{u(\tau)}\left(\theta,x\right)\right|dx\\\non
&&\le L\frac 1 {\theta}\int_{\tau}^{\tau+\theta}\liminf_{\eta\to 0}\frac{\|\bar U^{\sharp,\varepsilon}_{u(\tau)}\left(t-\tau+\eta\right)-S_\eta\bar U^{\sharp,\varepsilon}_{u(\tau)}\left(t-\tau\right)\|_{\L1(J_{t+\eta})}}{\eta}
\ee
The discontinuities of $\bar U^{\sharp,\varepsilon}_{u(\tau)}$ are solved by $S_\eta$ with exact shock or rarefaction for $x\not=0$ and with the $a$--Riemann solver in $x=0$ therefore the only difference between  $\bar U^{\sharp,\varepsilon}_{u(\tau)}\left(t-\tau+\eta\right)$ and $S_\eta\bar U^{\sharp,\varepsilon}_{u(\tau)}\left(t-\tau\right)$ are the rarefactions solved in an approximate way in the first function and in an exact way in the second. Recalling (\ref{eq:stima_rarefazioni}) we know that this error is of second order in the size of the rarefactions.\\
To show that (\ref{eq:first_int_condbappo}) holds, proceed as in (\ref{eq:49}). \square

%

%
%

\end{document}